\documentclass[a4wide]{article}

\usepackage[dvips]{epsfig} 
\usepackage[dvips]{hyperref}
\usepackage[latin1]{inputenc}
\usepackage{amsmath}
\usepackage{amsfonts}
\usepackage{amssymb}

\def\limfaible{\rightharpoonup} 

\title{${\cal T}$-class algorithms for pseudocontractions and $\kappa$-strict pseudocontractions
 in Hilbert spaces}
\author{Jean-Philippe Chancelier}
\begin{document}

\maketitle

\begin{abstract}In this paper we study iterative algorithms for finding a 
common element of the set of fixed points of $\kappa$-strict pseudocontractions 
or finding a solution of a variational inequality problem for a monotone, Lipschitz 
continuous mapping. The last problem being related to finding fixed points of 
pseudocontractions. These algorithms were already studied in \cite{acedo-xu} and 
\cite{nadezhkina} but our aim here is to provide the links between these know 
algorithms and the general framework of ${\cal T}$-class algorithms studied in \cite{bauschke-combettes}. 
\end{abstract}

%\begin{keyword}
% keywords here, in the form: keyword \sep keyword
%Nonexpansive mappings \sep Viscosity approximation \sep Fixed point.
% PACS codes here, in the form: \PACS code \sep code
%\end{keyword}

%\end{frontmatter}
\newtheorem{thm}{Theorem}
\newtheorem{cor}[thm]{Corollary}
\newtheorem{lem}[thm]{Lemma}
\newtheorem{claim}[thm]{Claim}
\newtheorem{axiom}[thm]{Axiom}
\newtheorem{conj}[thm]{Conjecture}
\newtheorem{fact}[thm]{Fact}
\newtheorem{hypo}[thm]{Hypothesis}
\newtheorem{assum}[thm]{Assumption}
\newtheorem{prop}[thm]{Proposition}
\newtheorem{crit}[thm]{Criterion}
\newtheorem{defn}[thm]{Definition}
\newtheorem{exmp}[thm]{Example}
\newtheorem{rem}[thm]{Remark}
\newtheorem{prob}[thm]{Problem}
\newtheorem{prin}[thm]{Principle}
\newtheorem{alg}{Algorithm}
\newtheorem{note}{Note}
\newtheorem{summ}{Summary}
\newtheorem{case}{Case}

\def\sequence#1{({#1}_n)_{n\ge 0}}
\def\Fix{\mathop{\normalfont Fix}}
\def\VI{\mathop{\normalfont VI}}
\def\defpar{\stackrel{\mbox{\tiny def}}{=}}
\def\argmax{\mathop{\mbox{\rm Argmax}}}
\def\argmin{\mathop{\mbox{\rm Argmin}}}
\def\bbP{{\mathbb P}} 
\def\bbC{{\mathbb C}} 
\def\bbE{{\mathbb E}} 
\def\E{{\cal E}}
\def\F{{\cal F}}
\def\H{{\cal H}}
\def\V{{\cal V}}
\def\W{{\cal W}}
\def\U{{\cal U}}
\def\R{{\mathbb R}}
\def\RB{{\mathbb R}}
\def\N{{\mathbb N}}
\def\M{{\mathbb M}}
\def\S{{\mathbb S}}
\def\bbR{{\mathbb R}} 
\def\bbI{{\mathbb I}} 
\def\bbU{{\mathbb U}} 
\def\Tad{{{\cal T}_{\mbox{\tiny ad}}}}
\def\texte#1{\quad\mbox{#1}\quad}
\def\Proba#1{\bbP\left\{ #1 \right\}} 
\def\Probax#1#2{{\bbP}_{#1}\left\{ #2 \right\}} 
\def\ProbaU#1#2{{\bbP}^{#1} \left\{ #2 \right\}} 
\def\ProbaxU#1#2#3{{\bbP}^{#1}_{#2} \left\{ #3 \right\}} 
\def\valmoy#1{\bbE\left[ #1 \right]}
\def\valmoyDebut#1{\bbE [ #1 } 
\def\valmoyFin#1{ #1 ]} 
\def\valmoyp#1#2{\bbE_{#1}\left[ #2 \right]}
\def\valmoypDebut#1#2{\bbE_{#1} \left[ #2 \right.} 
\def\valmoypFin#1{ \left. #1 \right]} 
\def\valmoypU#1#2#3{\bbE_{#1}^{#2}\left[ #3 \right]}
\def\norminf#1{ {\Vert #1 \Vert}_{\infty}}
\def\norm#1{ {\Vert #1 \Vert}}
\def\Hun{${\text{\bf H}}_1$}
\def\Hdeux{${\text{\bf H}}_2$}
\def\Htrois{${\text{\bf H}}_3$}
\def\psca#1{\left< #1 \right>}
\def\slim{\sigma\mbox{-}\lim}
%%\spnewtheorem*{proof}{Proof}{\itshape}{\rmfamily}

\newenvironment{myproof}{{\small{\it Proof~:}}}{\hfill$\Box$\normalsize
\\\smallskip}

\section{Introduction} 

Let $C$ be a closed convex subset of a Hilbert space $\H$ and  $P_C$ be the metric projection from 
$\H$ onto $C$. A mapping $Q : C \mapsto C$ is said to be a \emph{strict pseudocontraction} if there exists a
constant $0 \le \kappa < 1$ such that~:
\begin{equation}
 \norm{Q x - Q y}^2 \le \norm{x - y}^2 + \kappa \norm{(I - Q )x - (I - Q )y}^2\,,
 \label{spc}
\end{equation}
for all $x$, $y \in C$. A mapping $Q$ for which \eqref{spc} holds is also called 
a $\kappa$-strict pseudocontraction. As pointed out in \cite{acedo-xu} 
iterative methods for finding a common element of the set of fixed points of 
strict pseudocontractions are far less developed than iterative methods for 
nonexpansive mappings ($\kappa=0$) 
\cite{bib1,bib5,bib8,bib12,bib13,bib18,bib21,bib24,bib28,bib29,bib31}. 
We will, in section \ref{secacedo}
of this article, consider the algorithm \ref{firstalg} studied in \cite{acedo-xu} and we 
will show that this algorithm can be viewed as a ${\cal T}$-class algorithm as defined and studied 
in \cite{bauschke-combettes}. 

Section \ref{secnade} is devoted to the case $\kappa=1$ for which previous algorithm cannot be 
used. A mapping $A$ for which \eqref{spc} holds with $\kappa=1$ is called \emph{pseudocontractive}. 
We will see that \emph{pseudocontractive} mappings are related to monotone Lipschitz continuous
mappings. 
A mapping $A : C \mapsto \H$ is called \emph{monotone} if 
$$ \psca{Au - Av, u - v} \ge 0 \quad \mbox{for all} \quad (u, v) \in C^2\,.$$ 
$A$ is called $k$-Lipschitz continuous if there exists a positive real number $k$ such that
$$ \norm{Au - Av} \le k \norm{u - v} \quad \mbox{for all} \quad (u, v) \in C^2.$$
Let the mapping $A: C \mapsto \H$ be monotone and Lipschitz continuous.
The variational inequality problem is to find a $u \in  C$ such that 
$$ \psca {Au, v - u} \ge 0 \quad \mbox{for all} \quad  v \in C\,.$$
The set of solutions of the variational inequality problem is denoted by $\VI(C,A)$.

Assume that a mapping $Q : C \mapsto  C$ is pseudocontractive and $k$-Lipschitz-continuous then 
the mapping $A = I - Q$ is monotone and $(k + 1)$-Lipschitz-continuous and moreover 
$\Fix(Q) = \VI(C,A)$ \cite[Theorem 4.5]{nadezhkina} where $\Fix(Q)$ is the set of fixed points of $Q$, 
that is 
\begin{equation}
  \Fix(Q) \defpar  \left\{x \in C\,  :\,  Q x = x \right\}
\end{equation}

Thus, to cover the case $\kappa=1$, algorithms which aims at computing $P_{\VI(C,A)}x$ for 
a monotone and $k$-Lipschitz-continuous mapping $A$ are investigated. 
We will, in section \ref{secnade} mainly use results from  \cite{nadezhkina} to prove that the general 
algorithm that they use can be rephrased in a slightly extended ${\cal T}$-class algorithm framework.

\section{${\cal T}$-class iterative algorithm for a sequence of $\kappa$-strict pseudocontractions}

\label{secacedo}
Let $(Q_n)_{n\ge 0}$ be a sequence of $\kappa$-strict pseudocontractions, $\kappa \in [0,1)$ and 
$(\alpha_n)_{n\ge 0}$ a sequence of real numbers chosen so that $\alpha_n \in (\kappa,1)$. 
We consider as in \cite{acedo-xu} the following algorithm~:

\begin{alg}
\label{firstalg}
Given $x_0 \in C$, we consider the sequence $\sequence{x}$ generated by the following algorithm~:
\begin{eqnarray}
  &&y_n = \alpha_n x_n + (1 - \alpha_n) Q_n x_n, \nonumber \\ 
  &&C_n \defpar \left\{ z \in C \,| \, \norm{y_n - z}^2 \le \norm{x_n - z}^2 - (1 - \alpha_n)(\alpha_n - \kappa)
  \norm{x_n - Q_n x_n}^2 \right\}, \nonumber \\ 
  &&D_n \defpar \left\{ z \in C \,|\, \psca{x_n - z, x_0 - x_n} \ge 0 \right\}, \nonumber \\ 
  &&x_{n+1} = P_{(C_n \cap D_n)} x_0.\nonumber
\end{eqnarray} 
\end{alg}

We will show that this algorithm belong to the ${\cal T}$-class algorithms as defined in 
\cite{bauschke-combettes} and deduce its strong convergence to $P_F x_0$ when 
$F\ne \emptyset$ and where $F \defpar \cap_{n\ge 0} \Fix(Q_n)$. 

For $(x,y)\in \H^2$ define the mappings $H$ as follows~:
\begin{equation}
  H(x,y) \defpar \left\{ z \in \H \quad \vert \quad \psca{ z -y , x - y} \le 0  \right\}
\end{equation}
and denote by $Q(x,y,z)$ the projection of $x$ onto $H(x,y)\cap H(y,z)$. Note that $H(x,x)=\H$ and 
for $x \ne y$, $H(x,y)$ is a closed affine half space onto which $y$ is the projection of $x$. 

\begin{lem}The sequence generated by Algorithm \ref{firstalg} coincide with the sequence given 
by $x_{n+1}= Q(x_0,x_n,T_{n} x_n)$ with~: 
\begin{equation}
  \label{Tdef} 
  T_{n}(x) \defpar \frac{x+R_{n} y}{2} + 
  \frac{1}{2}\left(\frac{\kappa - \alpha_n}{1 -\alpha_n}\right)( x- R_{n} y), 
  \mbox{ and } R_{n}(x) \defpar \alpha_n x + (1-\alpha_n ) Q_n(x)\,.
\end{equation}
Moreover, we have~:
\begin{equation}
  2T_n-I = \kappa I +(1-\kappa) Q_nx \,.
  \label{Tdefn}
\end{equation}
\label{lemequiv}
\end{lem}

\begin{myproof}Let $\kappa \in [0,1)$,  $\alpha \in (\kappa,1)$, $y\defpar \alpha x + (1-\alpha ) Q x$ 
for a $\kappa$-strict pseudocontractions $Q$ and define $\Gamma(x,y)$ as follows~:
\begin{equation}
  \Gamma(x,y) \defpar \left\{ z \in \H \quad \, \vert \quad \, 
   \norm{y -z}^2 \le \norm{x -z }^2 -(1-\alpha)(\alpha -\kappa) \norm { x -Qx}^2 \right\} \,.
\end{equation}
We first prove that $\Gamma(x,y)= H(x,Tx)$ where $T$ is defined by equation \eqref{Tdef}. 
\begin{eqnarray}
&&  \norm{y -z }^2 - \norm{x -z }^2  \le  -(1-\alpha)(\alpha -\kappa) \norm { x -Qx}^2 \nonumber \\
&\Leftrightarrow&  \psca{y-z,y-z} - \norm{x -z }^2 \le  -(1-\alpha)(\alpha -\kappa) 
  \norm{x -Qx}^2 \nonumber \\
&\Leftrightarrow&  \psca{ y-x,y-z} + \psca{x-z,y- z} - \norm{x -z }^2  \le  -(1-\alpha)(\alpha -\kappa) 
  \norm{x -Qx}^2 \nonumber \\
&\Leftrightarrow & \psca{ y-x,y-z} + \psca{x-z,y- x}  \le  -(1-\alpha)(\alpha -\kappa) 
  \norm{x -Qx}^2 \nonumber \\
&\Leftrightarrow&  \psca{ y-x,y-z} + \psca{x-z,y- x}  \le  (\alpha -\kappa) \psca{ y-x , x-Qx} \nonumber \\
&\Leftrightarrow&  \psca{ y-x,y+x -2z + (\kappa -\alpha)( x-Qx) }  \le 0  \nonumber \\
&\Leftrightarrow&  \psca{ y-x,y+x -2z + \left(\frac{\kappa - \alpha}{1 -\alpha}\right){(x-y)}}  \le 0  \nonumber
\end{eqnarray}
which gives~: 
\begin{eqnarray}
  \psca{ z - \frac{x+y}{2} -\frac{1}{2} \left(\frac{\kappa - \alpha}{1 -\alpha}\right)(x-y), x-y} &\le& 0  \nonumber 
\end{eqnarray}
and since we have $x-Tx = (1/2)( 1 - \frac{\kappa - \alpha}{1 -\alpha})(x-y)$ with 
$( 1 - \frac{\kappa - \alpha}{1 -\alpha}) >0$ this is equivalent to $\psca{ z - Tx , x -Tx} \le 0$.
For $y_n =\alpha_n x_n + (1-\alpha_n ) Q_n x_n $, we thus obtain that $C_n = \Gamma(x_n,y_n) = H(x_n, T_n x_n)$ and 
since by definition of $H$ we have $D_n = H(x_0,x_n)$ the result follows. The last statement of the lemma \eqref{Tdefn} is obtained by simple rewrite from equation \eqref{Tdef}
\end{myproof}

We prove now that $T_n$ for all $n \in \N$ belongs to the ${\cal T}$ class of mappings. 
\begin{defn}
${\cal T} \defpar \left\{ T : \H \mapsto \H \, | \, \mbox{dom} T = \H \quad\mbox{and} \quad (\forall x \in \H)
  \Fix(T) \subset H(x,Tx)\right\}$ 
\end{defn}

\begin{lem}for all $n\in \N$ and $T_n$ defined by equation \eqref{Tdef} we have $T_n \in {\cal T}$.
\end{lem}

\begin{myproof}Using Lemma \ref{lemequiv} we have $2 T_n -I = \kappa I + (1-\kappa) Q_n$. If we can prove that 
when $Q$ is a $\kappa$-strict pseudocontraction the mapping $\kappa I + (1-\kappa) Q$ is quasi-nonexpansive 
then the result will follow from \cite[Proposition 2.3 (v)]{bauschke-combettes}. For $(x,y)\in \H^2$ we have~:
\begin{eqnarray}
  \Vert \kappa x &+& (1 -\kappa)Qx - y - (1-\kappa) y \Vert^2 = 
  \norm{ \kappa (x-y) + (1 -\kappa)(Qx -Qy) }^2 \nonumber \\ 
  &=& \kappa\norm{x-y}^2 + (1-\kappa)\norm{Qx -Qy}^2 - 
  \kappa(1-\kappa)\norm{x-y - (Qx-Qy)}^2 \nonumber \\  
  &=& \kappa\norm{x-y}^2 + (1-\kappa)\norm{Qx -Qy}^2 -  
  \kappa(1-\kappa)\norm{x-y - (Qx-Qy)}^2 \nonumber \\  
  &\le& \kappa\norm{x-y}^2 + (1-\kappa) \left( \norm{Qx -Qy}^2 -  
  \kappa\norm{(I-Q)x- (I-Q)y}^2 \right) \nonumber \\  
  &\le & \kappa\norm{x-y}^2 + (1-\kappa) \norm{x-y}^2 = \norm{x-y}^2\nonumber
\end{eqnarray}
Thus the mapping $\kappa I + (1-\kappa) Q$ is nonexpansive and thus also quasi-nonexpansive.
\end{myproof}

\begin{defn}\cite{bauschke-combettes} A sequence $\sequence{T}$ such that $T_n \in {\cal T}$ is {\emph coherent} 
if for every bounded sequence 
$\{z_n\}_{n\ge 0} \in \H$ there holds~: 
\begin{eqnarray}
  \left\{ 
  \begin{array}{l} 
    \sum_{n\ge 0} \norm{z_{n+1} -z_n}^2 < \infty \\
    \sum_{n\ge 0} \norm{z_{n} - T_n z_n}^2 < \infty 
  \end{array} 
  \right.
  \Rightarrow {\cal M}(z_n)_{n\ge 0} \subset \cap_{n \ge 0} \Fix(T_n)
  \label{coherentprop}
\end{eqnarray} 
\end{defn}
where ${\cal M}(z_n)_{n\ge 0}$ is the set of weak cluster points of the sequence $\sequence{z}$.

\begin{lem}Let $\sequence{Q}$ be a sequence of $\kappa$-strict pseudocontraction such that $\Fix(Q_n)=F$ which does not depends on $n$ and for each subsequence $\sigma{(n)}$ we can find a sub-sequence $\mu(n)$ such that 
$Q_{\mu(n)} \to Q$ with $\Fix(Q)=F$ and $Q$ is a $\kappa$-strict pseudocontraction. Then, the sequence $\sequence{T}$ given by \eqref{Tdef} is coherent. 
\end{lem}

\begin{myproof} Suppose that $\sequence{z}$ is a bounded sequence such that the left hand side of 
\eqref{coherentprop} is satisfied. 
Using \eqref{Tdefn} we have $\norm{z_{n} - T_n z_n}= (1 - \kappa )/2 \norm{z_n -Q_n z_n}$ and $\Fix(T_n)=\Fix(Q_n)$.
 Thus, verifying the coherence of $\sequence{T}$ or the coherence of $\sequence{Q}$ is equivalent. 
Consider now $u \in {\cal M}(z_n)_{n\ge 0}$, by hypothesis $\norm{z_{n} - Q_n z_n} \to 0$. Let $\sigma(n)$ a 
subsequence such that $z_{\sigma(n)} \limfaible u$, we extract a subsequence $\mu(n)$ such that $Q_{\mu(n)} \to Q $ 
and we thus obtain that $z_{\mu(n)} \limfaible u$ and $\norm{z_{\mu(n)} - Q z_{\mu(n)}} \to 0$. Now, if $Q$ 
is a $\kappa$-strict pseudocontraction, using \cite[Proposition 2.6]{acedo-xu} 
we have that $I-Q$ is demi-closed and thus $u \in \Fix(Q)=F$. 
\end{myproof}

\begin{rem} Given an integer $N \ge 1$, let, for each $1 \le i \le N$,  $S_i : C \mapsto C$ be a $\kappa_i$-strict 
pseudocontraction for some $0 \le \kappa_i < 1$. Let $\kappa \defpar \max\{\kappa_i \,:\, 1\le i \le N\}$. 
Assume the common fixed point set $ F \defpar \cap_{i=1}^N \Fix(S_i )$ of $\{S_i\}$ is nonempty. 
Assume also for each $n$, $\{\lambda_{n,i}\}_{i=1,\ldots,N}$ is a finite sequence of positive numbers 
such that $\sum_{i=1}^N \lambda_{n,i} = 1$ and $\inf_{n} \lambda_{n,i} >0$ for all $1 \le i \le N$. 
Let the mapping $Q_n : C \mapsto C $ be defined by~:
\begin{eqnarray}
  Q_n x \defpar \sum_{i=1}^N  \lambda_{n,i} S_i x \,.
\end{eqnarray}
Then using \cite{acedo-xu}, for all $n \in \N$, $Q_n$ is a $\kappa$-strict pseudocontraction and 
$\Fix(Q_n)=F$. Moreover for each subsequence $\lambda_{i,(\sigma_n)}$ we can extract a subsequence 
$\lambda_{i,\mu(n)}$ and $(\overline{\lambda_i})_{1 \le i\le N} \in (0,1)^N$ such that 
$\lambda_{i,\mu(n)} \to \overline{\lambda_i}$ for all $1 \le i\le N$. We thus have 
$Q_{\mu(n)} \to \sum_{i} \overline{\lambda_i} S_i $ and using previous lemma 
the sequence $\sequence{T}$ is coherent.
\label{acedoxu}
\end{rem}

Given $T_n \in {\cal T}$ we can also consider \cite{bauschke-combettes} the following algorithm~:
\begin{alg}Given $\epsilon \in (0,1]$ and $x_0 \in C$ we consider the sequence given by the iterations 
$x_{n+1} = x_n + (2-\epsilon)(T_n x_n - x_n)$. 
\label{secondalg}
\end{alg}
Gathering previous result the strong convergence of Algorithm \ref{firstalg} to $P_F x_0$ and the weak 
convergence of Algorithm \ref{secondalg} is obtained by 
\cite[Theorem 4.2]{bauschke-combettes} that we recall now~:

\begin{thm}\cite[Theorem 4.2]{bauschke-combettes} Suppose that $(T_n)_{n\ge0}$ is coherent. Then  \\
\noindent $(i)$ if $F\ne \emptyset$, then every orbit of Algorithm \ref{secondalg} converges weakly to a point in $F$ \\
\noindent$(ii)$ For an arbitrary orbit of Algorithm \ref{firstalg}, exactly one of the following alternatives holds~: 
\begin{enumerate}
  \item[(a)] $F\ne \emptyset$ and $x_n \to_{n} P_F x_0$.
  \item[(b)] $F=  \emptyset$ and $x_n \to_{n} +\infty $.
  \item[(c)] $F=  \emptyset$ and the algorithm terminates. 
\end{enumerate}
\end{thm}

\begin{rem}Note that using previous theorem and Remark \ref{acedoxu} we obtain an other proof of 
\cite[Theorem 5.1]{acedo-xu}. In fact the proofs are very similar but we just hilite here the role 
played by ${\cal T}$-class sequences. 
\end{rem}

\section{${\cal T}$-class iterative algorithm for a sequence of pseudo contractions}

\label{secnade}

Let $F$ be a closed convex of $\H$ we define ${\cal U}_F$ as follows~: 
\begin{equation} 
  {\cal U}_F \defpar \left\{ T : \H \mapsto \H \, | \, \mbox{dom} T = \H 
  \quad\mbox{and} \quad (\forall x \in \H)
  F \subset H(x,Tx)\right\}\,.
\end{equation} 
Of course we have $T\in {\cal T}\Leftrightarrow T \in {\cal U}_{\Fix(T)}$. 

A mapping $Q : \H \mapsto \H$ is said $F$-quasi-nonexpansive if 
\begin{equation}
  \forall (x,y) \in \H \times F \quad \norm{Q x -y } \le \norm{x -y} 
\end{equation}
and we can characterize elements of ${\cal U}_F$ using the following easy lemma~:
\begin{lem} $2T -I$ is $F$-quasi-nonexpansive is equivalent to $T \in {\cal U}_F$. 
\label{lemfquasi}
\end{lem}
\begin{myproof}The proof follows from the equality \cite[(2.6)]{bauschke-combettes}~: 
\begin{equation} 
  (\forall (x,y) \in \H^2) \quad 4 \psca{y -Tx, x- Tx} = \norm{ (2T-I)x -y}^2 - \norm{x-y}^2\,.
\end{equation} 
\end{myproof}

\begin{defn}A sequence $\{T_n\}_{n\ge 0} \subset {\cal U}_F$ is $F$-{\emph coherent} if 
for every bounded sequence $\{z_n\}_{n\ge 0} \in \H$ there holds~: 
\begin{eqnarray}
  \left\{ 
  \begin{array}{l} 
    \sum_{n\ge 0} \norm{z_{n+1} -z_n}^2 < \infty \\
    \sum_{n\ge 0} \norm{z_{n} - T_n z_n}^2 < \infty 
  \end{array} 
  \right.
  \Rightarrow {\cal M}(z_n)_{n\ge 0} \subset F 
  \label{coherentprop1}
\end{eqnarray} 
\end{defn}

We propose now the following extension of \cite[Theorem 4.2]{bauschke-combettes} for the two 
algorithms \ref{secondalg} and \ref{firstalg-comb}.

\begin{alg}Given $x_0 \in C$ we consider the sequence given by the iterations 
  $$x_{n+1}= Q(x_0,x_n,T_{n} x_n)$$
\label{firstalg-comb} 
\end{alg}

\begin{thm}Suppose that $(T_n)_{n\ge0}$ is $F$-coherent for a closed convex $F$ Then 
$(i)$ if $F\ne \emptyset$, then every orbit of Algorithm \ref{secondalg} converges 
weakly to a point in $F$ $(ii)$ For an arbitrary orbit of Algorithm \ref{firstalg-comb}, 
exactly one of the following alternatives holds~: 
\begin{enumerate}
  \item[(a)] $F\ne \emptyset$ and $x_n \to_{n} P_F x_0$.
  \item[(b)] $F=  \emptyset$ and $x_n \to_{n} +\infty $.
  \item[(c)] $F=  \emptyset$ and the algorithm terminates. 
\end{enumerate}
\label{fcoherentthm}
\end{thm}

\begin{myproof}The result is very similar to \cite[Theorem 2.9]{bauschke-combettes} and a 
careful reading of the proof and remarks in \cite{bauschke-combettes,combettes} leads to 
the conclusion that it remains true as stated here. 
\end{myproof}

We give now a typical application of this theorem. 

\def\Tu{T^{(1)}}
\def\Td{T^{(2)}}
\def\Tul{\Tu_{\lambda}}
\def\Tdl{\Td_{\lambda}}
\def\Tun{\Tu_{n}}
\def\Tdn{\Td_{n}}
\begin{defn}
\label{Top}
For $A: C \mapsto C$ a monotone and $k$-Lipschitz mapping, let $T_{\lambda} : \H \times \H \mapsto \H$ the 
mapping defined by $T_{\lambda}(x,y) \defpar P_C(x -\lambda Ay)$. 
We also define $\Tul x \defpar T_{\lambda}(x,x)$ and $\Tdl x \defpar T_{\lambda}(x,T_{\lambda}(x,x))= T_{\lambda}(x,\Tul x)$. 
\label{Tddef}
\end{defn}

We assume that $\lambda k \in [a,b] \subset (0,1)$ and consider $(\lambda_n)_{n \ge 0}$ a sequence of 
real numbers such that $\lambda_n k \in [a,b]$. To simplify the notations 
we will use $\Tun$ (resp. $\Tdn$) for denoting $\Tu_{\lambda_n}$ (resp. $\Td_{\lambda_n}$). 

Let $F \defpar \VI(C,A)$, It is known that $F$ is closed convex and that we have $\Fix{\Tul}= F$. 
It is easy to see that $F \subset \Fix(\Tdl)$ but the inclusion may be strict and thus we do not expect 
the mapping $\Tdl$ to be quasi-nonexpansive. Following inequalities contained in the proof of 
\cite[Theorem 3.1]{nadezhkina} we obtain $F$-quasi-nonexpansive property as exposed now.

\begin{lem}$\Tdl$ is $F$-quasi-nonexpansive where $F \defpar \VI(C,A)$ or using 
  Lemma \ref{lemfquasi} $(\Tdl +I)/2 \in {\cal U}_F$.
  \label{lemquasinonexp}
\end{lem}

\begin{myproof} Let $y =\Tul(x)$ and $u \in \VI(C,A)$. We use the fact that for all $x\in \H$ and $y\in C$ 
$P_C x$ can be characterized as follows~:
\begin{equation}
  \norm{x-y}^2 \ge   \norm{x- P_Cx}^2  + \norm{y -P_C x}^2 
\end{equation}
and since $A$ is a monotone mapping following the steps of the proof of \cite[Theorem 3.1]{nadezhkina} 
that we reproduce here we obtain~:
\begin{eqnarray} 
  \norm{\Tdl(x) - u}^2 &\le&  \norm{x - \lambda A y - u}^2 - \norm{x - \lambda A y - \Tdl(x)}^2 \nonumber \\
  &=& \norm{x - u}^2 - \norm{x - \Tdl(x)}^2 + 2\lambda \psca{A y, u - \Tdl(x)} \nonumber \\
  &=& \norm{x - u}^2 - \norm{x - \Tdl(x)}^2  \nonumber \\
  & & \quad + 2\lambda ( \psca{A y - Au, u - y} + \psca{Au, u - y} + \psca{A y, y - \Tdl(x)})\nonumber \\
  &\le& \norm{x - u}^2 - \norm{x - \Tdl(x)}^2 + 2\lambda \psca{A y, y - \Tdl(x)} \nonumber \\
  &=& \norm{x - u}^2 - \norm{x - y}^2 - 2 \psca{x - y, y - \Tdl(x)} - \norm{y - \Tdl(x)}^2\nonumber \\
  & & \quad +2\lambda \psca{A y, y - \Tdl(x)} \nonumber \\
  &=& \norm{x - u}^2 - \norm{x - y}^2 - \norm{y - \Tdl(x)}^2  \nonumber \\ 
  && + 2\psca{x - \lambda A y - y, \Tdl(x) - y}\, . \nonumber 
\end{eqnarray} 
Further, since $y = P_C (x - \lambda A x)$ and $A$ is $k$-Lipschitz-continuous, we have
\begin{eqnarray} 
  \left< x -  \lambda A y - y \right. &,& \left. \Tdl(x) - y \right> = \psca{x - \lambda A x - y, \Tdl(x) - y} \nonumber \\
    &+& \psca{\lambda A x - \lambda A y, \Tdl(x) - y} \le \psca{\lambda Ax - \lambda A y, \Tdl(x) - y} \nonumber \\
  &\le& \lambda k \norm{x - y} \norm{\Tdl(x) - y}\,.\nonumber
\end{eqnarray} 
So, we have~;
\begin{eqnarray} 
  \norm{\Tdl(x) - u}^2 &\le&  \norm{x - u}^2 - \norm{x - y}^2 - \norm{y - \Tdl(x)}^2 + 2\lambda k \norm{x - y}\norm{\Tdl(x) - y}\nonumber \\
  &\le& \norm{x - u}^2 + (\lambda^2 k^2 - 1) \max\left( \norm{x - y}^2, \norm{\Tdl(x) - y}^2 \right) \label{equay} \\
  &\le& \norm{x - u}^2 \,.\nonumber 
\end{eqnarray} 
\end{myproof}

\begin{cor}If we consider $R \defpar \alpha I + (1-\alpha)S$ where $S$ is a non-expansive mapping and 
define $\tilde{F} = \Fix(S)\cap \VI(C,A)$ then we obtain immediately that $R\circ \Tdl$ is a 
$\tilde{F}$-quasi-nonexpansive mapping. 
\end{cor}

\begin{myproof}Let $u \in \tilde{F}$ then $u=Ru$ and we have $\norm{R\circ \Tdl - u} \le \norm{\Tdl -u} $ and 
the previous lemma ends the proof.
\end{myproof}

\begin{lem}The sequence $Q_n = 1/2(\Td_n+I)$ is $F$-coherent. 
\end{lem}

\begin{myproof}Let $(y_n)_{n\ge 0}$ a bounded sequence satisfying the left hand side of equation 
  \eqref{coherentprop1} and $\varphi \in {\cal M}(y_n)_{n\ge 0}$. We can find a subsequence $y_{\sigma(n)}$ which converges 
  weakly to $\varphi$. For simplicity, we use the notation $y_n$ for the subsequence and since it satisfies the 
  left hand side of equation \eqref{coherentprop1} we have $\norm{y_n - Q_n y_n}\to 0$. By definition of $Q_n$ 
  we also have $\norm{y_n - \Td_n y_n}\to 0$ and thus $\Td_n y_n \rightharpoonup u$
  From equation \eqref{equay}
  we obtain~: 
  \begin{eqnarray} 
  \norm{\Tdl x - u}^2 \le  \norm{x - u}^2 + (\lambda^2 k^2 - 1)  \max\left( \norm{x - \Tul x}^2, \norm{\Tdl x - \Tul x}^2 \right)\label{ineqmax1} \nonumber
\end{eqnarray} 
Thus~:
\begin{eqnarray} 
 \max\left( \norm{x - \Tul x}^2 \right. &,& \left. \norm{\Tdl x - \Tul x}^2 \right) \le  \frac{1}{1 - \lambda^2 k^2} \left( \norm{x - u}^2 -  \norm{\Tdl x - u}^2 \right) \nonumber \\
 & \le & K \left( \norm{x - u} + \norm{\Tdl x -u} \right) \norm{ x - \Tdl x} \label{ineqmax}
\end{eqnarray} 
Using Lemma \ref{lemquasinonexp}, the sequence $ \Td_n y_n$ is bounded and we thus have from the previous inequality 
$\norm{y_n - \Tu_n y_n} \to 0$ and $\norm{\Td_n y_n - \Tu_n y_n} \to 0$. 

Using next lemma (Lemma \ref{lemvi}) we therefore obtain that for $(v,w)\in G(T)$~: 
$$ \psca{v - \varphi , w} =  \lim_{n \to \infty} \psca{v - \Td_n y_n , w} \ge 0\,. $$
Thus we obtain that $\psca{v - \varphi , w} \ge 0$  which gives $\varphi \in T^{-1}(0)$ since $T$ is maximal monotone and then $\varphi \in F=\VI(C,A)$. Thus
$Q_n$ is $F$-coherent.
\end{myproof}

\begin{cor}Let $(R_n)_{n\ge 0}$ a sequence of nonexpansive mappings such that for each subsequence 
$\sigma(n)$ it is possible to extract a subsequence $\mu(n)$ and find $R_\mu$ such that 
$R_{\mu(n)} y_n \to_{n\to \infty} R_{\mu} y_n$ for every bounded sequence $(y_n)_{n\ge 0}$ 
with $\Fix{R_{\mu}}= {\cal S}$ a fixed set such that ${\cal S} \cap {\cal S} \ne \emptyset$. Then, 
we also have that $Q_n = 1/2((R_n\circ \Td_n)+I)$ is $F\cap {\cal S}$-coherent. 
\label{Rcor}
\end{cor}

\begin{myproof}Let $u \in {\cal S} \cap {\cal S}$, since $R_n$ is nonexpansive we have~: 
$\norm{R_n \circ \Tdl - u} \le \norm{\Tdl -u}$, Thus equation \eqref{ineqmax} can be replaced by~: 
  \begin{eqnarray} 
    \norm{R_n \circ \Tdl x - u}^2 \le  \norm{x - u}^2 + (\lambda^2 k^2 - 1)  \max\left( \norm{x - \Tul x}^2, \norm{\Tdl x - \Tul x}^2 \right) \nonumber
\end{eqnarray} 
proceeding as in previous lemma we obtain that for $(y_n)_{n\ge 0}$ a bounded sequence satisfying the left hand side of equation \eqref{coherentprop1} for the sequence of mapping $R_n \circ \Td_n$ we also have up to 
subsequences that $\norm{y_n - \Tu_n y_n} \to 0$ and $\norm{\Td_n y_n - \Tu_n y_n} \to 0$ and thus also 
$\norm{y_n - \Td_n y_n} \to 0$. Thus, as before,
 if $\varphi$ is a weak limit of $(y_n)_{n\ge 0}$ we have $\varphi \in F$. 
Moreover, we have~: 
\begin{eqnarray} 
  \norm{\Td_n y_n - R_{\mu} \nu } 
  &\le&  
  \norm{\Td_n y_n - y_n} + \norm{y_n - R_{n}\circ\Td_n y_n} \nonumber \\
  &&+ \norm{ R_{n} \circ \Td_n y_n - R_{\mu} \circ \Td_n y_n }
  + \norm{ \Td_n y_n - \nu }
\end{eqnarray} 
Thus~
$$
\liminf_{n\mapsto \infty} \norm{\Td_n y_n - R_{\mu} \nu } \le  \liminf_{n\mapsto \infty} 
\norm{ \Td_n y_n - \nu }$$
which by Opial's condition is only possible if $R_{\mu} \nu = \nu$. We conclude that 
$\nu \in F\cap {\cal S}$ which ends the proof.
\end{myproof}

\begin{lem}\cite{nadezhkina} Let $T : \H \mapsto H$ the mapping defined by $Tv \defpar Av + N_Cv$ when 
$v \in C$ and $Tv = 0$ when $v \not\in C$ where $N_C$ is the normal cone to $C$ at $v \in C$. 
Let $G(T)$ be the graph of $T$ and $(v,w) \in G(T)$. Then for $x \in C$ we have the following 
inequality~: 
\begin{eqnarray}
  \psca{v - \Tdl x , w} \ge \psca{v -\Tdl x,A \Tdl x - A \Tul x} - 
  \psca{v - \Tdl x ,\frac{\Tdl x - x}{\lambda}} \nonumber
\end{eqnarray} 
\label{lemvi}
\end{lem}
\begin{myproof}The proof of this inequality is given in \cite{nadezhkina}, we reproduce it for 
the sake of completeness. The mapping $T$ is maximal monotone, and $0 \in Tv$  if and only if 
$v \in  \VI(C,A)$. Let $(v,w) \in G(T)$. Then, we have $w \in Tv = Av + N_Cv$ and hence 
$w - Av \in N_Cv$. So, we have $ \psca{v - t,w - Av} \ge 0$ for all
$t \in C$ . On the other hand, from $\Tdl(x)= P_C (x - \lambda A\Tul(x))$ and $v \in  C$ we have
$\psca{x - \lambda A y -\Tdl(x), \Tdl(x) - v} \ge 0$ and hence $\psca{v - \Tdl(x) ,\Tdl(x) - x \lambda + A\Tul x} \ge 0$.
From $\psca{v - t,w - Av} \ge 0$ for all $t \in C$ and $\Tdl(x) \in C$, we have
\begin{eqnarray}
  \psca{v - \Tdl x, w} &\ge& \psca{v - \Tdl x,Av}\nonumber \\
  &\ge& \psca{v - \Tdl x,Av} - \psca{v - \Tdl x , \frac{ \Tdl x - x}{\lambda} + A \Tul x} \nonumber \\
  & = &  \psca{v - \Tdl x,Av - A \Tdl x } + \psca{v - \Tdl x,A\Tdl x - A \Tul x}\nonumber \\
  && - \norm{v - \Tdl x , \frac{\Tdl x - x}{\lambda}}\nonumber \\
  & \ge & \psca{v -\Tdl x,A \Tdl x - A \Tul x} - \psca{v - \Tdl x , \frac{\Tdl x - x}{\lambda}} \nonumber
\end{eqnarray} 
\end{myproof}

We end this section by gathering previous results in a main theorem. 
The proof is immediate by applying Theorem \ref{fcoherentthm}. The first statement is a new 
result. The second statement when applied to the sequence $R_n = \alpha_n Id + (1 - \alpha_n) S$ with 
$\alpha_n \in [0,c)$ and $c < 1$ gives the same result as \cite[Theorem 3.1]{nadezhkina}. 

\begin{thm}Let $(R_n)_{n\ge 0}$ a sequence of nonexpansive mappings satisfying the 
hypothesis of Corollary \ref{Rcor} and $(\Td_n)_{n\ge 0}$ the sequence of mappings 
defined on Definition \ref{Tddef}. Then, every orbit of Algorithm \ref{secondalg} 
applied to the sequence of mappings $R_n \circ \Td_n$ converges weakly to a point in $F$ 
and the sequence generated by Algorithm \ref{firstalg} 
converges strongly to $P_F x_0$. 
\end{thm}

%\bibliographystyle{elsart-num-sort.bst}
%\bibliography{non-exp-map-fejer-jpc}

\begin{thebibliography}{10}
\expandafter\ifx\csname url\endcsname\relax
  \def\url#1{\texttt{#1}}\fi
\expandafter\ifx\csname urlprefix\endcsname\relax\def\urlprefix{URL }\fi

\bibitem{acedo-xu}
G.~L. Acedo, H.-K. Xu, Iterative methods for strict pseudo-contractions in
  hilbert spaces, Nonlinear Analysis 67 (2007) 2258--2271.

\bibitem{bib1}
H.~Bauschke, The approximation of fixed points of compositions of nonexpansive
  mappings in hilbert space, J. Math. Anal. Appl. 202 (1996) 150--159.

\bibitem{bauschke-combettes}
H.~H. Bauschke, P.~L. Combettes, A weak-to-strong convergence principle for
  fejér-monotone methods in hilbert spaces, Mathematics of Operations Research
  26~(2) (2001) 248--264.

\bibitem{combettes}
P.~Combettes, S.~Histoaga, Equilibrium programming in hilbert spaces, Journal
  of Nonlinear and Convex Analysis 6~(1) (2005) 117--136.

\bibitem{bib5}
K.~Goebel, W.~Kirk, Topics in Metric Fixed Point Theory, Cambridge Studies in
  Advanced Mathematics, {Cambridge University Press} ed., 1990.

\bibitem{bib8}
B.~Halpern, Fixed points of nonexpanding maps, Bull. Amer. Math. Soc. 73 (1967)
  957--961.

\bibitem{bib12}
T.~Kim, H.~Xu, Strong convergence of modified mann iterations for
  asymptotically nonexpansive mappings and semigroups, Nonlinear Anal. 64
  (2006) 1140--1152.

\bibitem{bib13}
P.~Lions, Approximation de points fixes de contractions, C. R. Acad. Sci. Série
  {A--B} Paris 284 (1977) 1357--1359.

\bibitem{nadezhkina}
N.~Nadezhkina, W.~Takahashi, strong convergence theorem by a hybrid method for
  nonexpansive mappings and lipschitz-continuous monotone mappings, siam j.
  optim 16~(4) (2006) 1230--1241.

\bibitem{bib18}
K.~Nakajo, W.~Takahashi, Strong convergence theorems for nonexpansive mappings
  and nonexpansive semigroups, J. Math. Anal. Appl. 279 (2003) 372--379.

\bibitem{bib21}
S.~Reich, Weak convergence theorems for nonexpansive mappings in banach spaces,
  J. Math. Anal. Appl. 67 (1979) 274--276.

\bibitem{bib24}
N.~Shioji, W.~Takahashi, Strong convergence of approximated sequences for
  nonexpansive mappings in banach spaces, Proc. Amer. Math. Soc. 125 (1997)
  3641--3645.

\bibitem{bib28}
R.~Wittmann, Approximation of fixed points of nonexpansive mappings, Arch.
  Math. 58 (1992) 486--491.

\bibitem{bib29}
H.~Xu, Iterative algorithms for nonlinear operators, J. London Math. Soc. 66
  (2002) 240--256.

\bibitem{bib31}
H.~Xu, Strong convergence of an iterative method for nonexpansive mappings and
  accretive operators, J. Math. Anal. Appl. 314 (2006) 631--643.

\end{thebibliography}

\end{document}